\documentclass[a4paper,12pt,leqno]{amsart}
\usepackage{amsmath,amssymb,amsxtra,amsfonts,comment,graphicx,psfrag}
\usepackage{bm,mathrsfs}
\usepackage{mathtools}
\usepackage{xspace}
\usepackage{color}
\usepackage{epsfig}

\usepackage{enumerate}
\usepackage{hyperref}
\newdimen\tempdima
\newdimen\tempdimb
\newbox\tempboxa
\newdimen\GGGlength
\newdimen\GGGheight
\newdimen\GGGborder
\newbox\GGGbox

\def\GGGput[#1,#2](#3,#4)#5{%
  \setbox\GGGbox\vbox{\hbox{#5}\kern0pt}%
  \GGGlength\wd\GGGbox%
  \divide\GGGlength by100 \multiply\GGGlength by#1%
  \GGGheight\ht\GGGbox%
  \divide\GGGheight by100 \multiply\GGGheight by#2%
  \put(#3,#4){\kern-\GGGlength\raise-\GGGheight\box\GGGbox}}

\def\GGGputRG[#1,#2](#3)(#4,#5)[#6]#7{%
  \setbox\GGGbox\vbox{\hbox{#7}\kern0pt}%
  \GGGlength\wd\GGGbox\tempdima\GGGlength%
  \divide\GGGlength by100 \multiply\GGGlength by#1%
  \GGGheight\ht\GGGbox\tempdimb\GGGheight%
  \divide\GGGheight by-100 \multiply\GGGheight by#2%
  \advance\GGGheight\tempdimb%
  \advance\tempdima\GGGborder\advance\tempdima\GGGborder%
  \advance\tempdimb\GGGborder\advance\tempdimb\GGGborder%
  \setbox\tempboxa\hbox{\vbox to 0mm{\vss\hbox to 0mm{\hss%
\graysqr{\number\tempdima}{\number\tempdimb}{#3}\kern\GGGborder\hfill}\kern-\GGGborder}%
    \box\GGGbox}%
  \setbox\tempboxa\vbox to
0cm{\kern-\GGGheight\hbox{\kern-\GGGlength\box\tempboxa}\vss}%
  \put(#4,#5){\rotatebox{#6}{\box\tempboxa}}%
  }

\def\graysqr#1#2#3{\immediate{
}}

\def\iu{{\rm i}}
\def\e{{\rm e}}
\def\d{{\rm d}}
\def\Z{{\mathbb Z}}

\def\R{{\mathbb R}}
\def\C{{\mathbb C}}
\def\integer{{\mathbb Z}}
\def\real{{\mathbb R}}
\def\torus{{\mathbb T}}
\def\eps{\varepsilon}
\def\bigo{{\mathcal O}}
\newcommand{\sfrac}[2]{\mbox{\footnotesize$\displaystyle\frac{#1}{#2}$}}
\def\M{\mathcal{M}}
\newcommand\calM{\mathcal{M}}
\newcommand\calH{\mathcal{H}}
\newcommand\calL{\mathcal{L}}
\newcommand\calE{\mathcal{E}}
\newcommand\calU{\mathcal{U}}

\newcommand\bfk{k}

\newcommand\bfy{{\mathbf y}}
\newcommand\bfz{{\mathbf z}}

\newcommand\bfomega{\omega}

\newcommand\bfzero{0}

\def\for{\quad\hbox{ for}\quad }
\def\with{\quad\hbox{ with}\quad }
\def\eps{\varepsilon}

\def\phi{\varphi}

\def\dt{\Delta t}
\def\dx{\Delta x}

\DeclareMathOperator\Imag{Im}

\theoremstyle{plain}

\theoremstyle{definition}

\numberwithin{equation}{section}
\numberwithin{table}{section}
\numberwithin{figure}{section}

\begin{document}

\title[Dynamics, numerical analysis, and some geometry]
{Dynamics, numerical analysis, \\ and some geometry}

\author[Ludwig Gauckler]{Ludwig Gauckler}
\address{Institut f\"ur Mathematik, FU Berlin,
Arnimalle 9,  D-14195 Berlin, Germany.} 
\email {\href{mailto:gauckler@math.fu-berlin.de}{gauckler{\it @\,}math.fu-berlin.de}} 

\author[Ernst Hairer]{Ernst Hairer}
\address{Section de math\'ematiques, 2-4 rue du Li\`evre,
Universit\'e de Gen\`eve, CH-1211 Gen\`eve 4, Switzerland.} 
\email {\href{mailto:Ernst.Hairer@unige.ch}{Ernst.Hairer{\it @\,}unige.ch}} 

\author[Christian Lubich]{Christian Lubich}
\address{Mathematisches Institut, Universit\"at T\"ubingen, Auf der Morgenstelle 10, 
D-72076 T\"ubingen, Germany.}
\email {\href{mailto:lubich@na.uni-tuebingen.de}{lubich{\it @\,}na.uni-tuebingen.de}}  


\date{October 10, 2017}

\begin{abstract} Geometric aspects play an important role in the construction and analysis of structure-preserving numerical methods for a wide variety of ordinary and partial differential equations. Here we review the development and theory of symplectic integrators for Hamiltonian ordinary and partial differential equations, of dynamical low-rank approximation of time-dependent large matrices and tensors, and its use in numerical integrators for Hamiltonian tensor network approximations in quantum dynamics.
\end{abstract}

\maketitle

\section{Introduction} 

It has become a commonplace notion in all of numerical analysis (which here is understood as comprising the construction and the mathematical analysis of numerical algorithms) that a good algorithm should ``respect the structure of the problem'' --- and in many cases the ``structure'' is of {\it geometric} nature. This immediately leads to two basic questions, which need to be answered specifically for each problem:
\begin{itemize}
\item How can  numerical methods be constructed that ``respect the geometry'' of the problem at hand?
\item What are  benefits from using a structure-preserving algorithm for this problem, and how do they come about?
\end{itemize}
In this note we present results in the numerical analysis  of {\it dynamic} (evolutionary, time-dependent) ordinary and partial differential equations for which geometric aspects play an important role. These results belong to the area that has become known as {\it Geometric Numerical Integration}, which has developed vividly in the past quarter-century, with substantial contributions by researchers with very different mathematical backgrounds. We just refer  to the books (in chronological order)
\cite{sanz-serna94nhp,hairer02gni,suris03tpo,leimkuhler04shd,hairer06gni,lubich08fqt,feng10sga,faou12gni,wu13spa,blanes16aci}  
and to the Acta Numerica review articles \cite{sanz-serna92sif,iserles00lm,marsden01dma,mclachlan02sm,hairer03gni,deckelnick05cog,bond07mda,chu08laa,hochbruck10ei,wanner10kna,christiansen11tis,abdulle12thm,dziuk13fem}. In this note we restrict ourselves to some selected topics to which we have contributed. 

In Section~\ref{sec:ham-ode} we begin with reviewing numerical methods for approximately solving {\it Hamiltonian systems} of ordinary differential equations, which are ubiquitous in many areas of physics. Such systems are characterized by the {\it symplecticity} of the flow, a geometric property that one would like to transfer to the numerical discretization, which is then called a {\it symplectic integrator}. Here, the two questions above become the following:
\pagebreak[3]
\begin{itemize}
\item How are symplectic integrators constructed?
\item What are favourable long-time properties of symplectic integrators, and how can they be explained?
\end{itemize}
The first question relates numerical methods with the theories of Hamilton and Jacobi from the mid-19th century, and the latter question connects numerical methods with the analytical techniques of Hamiltonian perturbation theory, a subject developed from the late 19th throughout the 20th century, from Lindstedt and Poincar\'e and Birkhoff to Siegel and Kolmogorov, Arnold and Moser (KAM theory), to Nekhoroshev and further eminent mathematicians. This connection comes about via {\it backward error analysis}, which is a concept that first appeared  in numerical linear algebra \cite{wilkinson60eao}.  The viewpoint is to interpret the numerical approximation as the exact (or almost exact) solution of a modified equation. In the case of a symplectic integrator applied to a Hamiltonian differential equation, the modified differential equation turns out to be again Hamiltonian, with a Hamiltonian that is a small perturbation to the original one. This brings Hamiltonian perturbation theory into play for the long-time analysis of symplectic integrators.
Beyond the purely mathematical aspects, it should be kept in mind that symplectic integrators are first and foremost an important tool in computational physics. In fact such numerical methods appeared first in the physics literature \cite{vogelaere56moi,verlet67ceo,ruth83aci}, in such areas as nuclear physics and molecular dynamics, and slightly later  \cite{wisdom91smf} in celestial mechanics, which has been the original motivation in the development of Hamiltonian perturbation theory \cite{poincare92lmn,siegel71vuh}. It was not least with the use of symplectic integrators that the centuries-old question about the stability of the solar system was finally answered negatively in the last decade by Laskar; see \cite{laskar13its} and compare also with \cite{moser78its}.

In Section~\ref{sec:ham-osc} we consider numerical methods for finite-dimensional {\it Hamiltonian systems with multiple time scales}  where, in the words of Fermi, Pasta \& Ulam \cite{fermi55},  ``the non-linearity is introduced as a perturbation to a primarily linear problem. The behavior of the systems is to be studied for times which are long compared to the characteristic periods of the corresponding linear problem.'' The two basic questions above are reconsidered for such systems. Except for unrealistically small time steps, the backward error analysis of Section~\ref{sec:ham-ode} does not work for such systems, and a different technique of analysis is required.   {\it Modulated Fourier expansions} in time were originally developed (since 2000) for studying numerical methods for such systems and were subsequently  also recognized as a powerful analytical technique for proving new results  for continuous systems of this type, including the original Fermi--Pasta--Ulam system. While the canonical transformations of Hamiltonian perturbation theory transform the system into a normal form from which long-time behaviour can be read off, modulated Fourier expansions embed the system into a high-dimensional system that has a Lagrangian  structure with invariance properties that enable us to infer long-time properties of the original system. Modulated Fourier expansions do not use nonlinear coordinate transformations, which is one reason for their suitability for studying numerical methods, which are most often not invariant under nonlinear transformations.

In Section~\ref{sec:ham-pde} we present long-time results for suitable numerical discretizations of {\it Hamiltonian partial differential equations} such as nonlinear wave equations and nonlinear Schr\"odinger equations. A number of important results on this topic have been obtained in the last decade, linking the numerical analysis of such equations to recent advances in their mathematical analysis. 
The viewpoint we take here is to consider the Hamiltonian partial differential equation as an infinite-dimensional system of the oscillatory type of Section~\ref{sec:ham-osc} with infinitely many frequencies, and we present results on the long-time behaviour of the numerical and the exact solutions that have been  obtained with modulated Fourier expansions or with techniques from infinite-dimensional Hamiltonian perturbation theory. We mention, however, that there exist other viewpoints on the equations considered, with different geometric concepts such as multisymplecticity \cite{bridges97msa,marsden98mgv}. While  multisymplectic integrators, which preserve this geometric structure, have been constructed and favourably tested in numerical experiments \cite{bridges01msi,ascher04mbs} (and many works thereafter), as of now there appear to be no proven results on the long-time behaviour of such methods.

In Section~\ref{sec:dlr} we consider {\it dynamical low-rank approximation}, which leads to a different class of dynamical problems with different geometric aspects. The problem here is to approximate large (or rather too large, huge) time-dependent matrices, which may be given explicitly or are the unknown solution to a matrix differential equation, by matrices of a prescribed rank, typically much smaller than the matrix dimension so that a data-compressed approximation is obtained. Such problems of data and/or model reduction arise in a wide variety of applications ranging from information retrieval to quantum dynamics. 
On projecting the time derivative of the matrices to the tangent space of the manifold of low-rank matrices at the current approximation, this problem leads to a differential equation on the low-rank manifold, which needs to be solved numerically.  We present answers to the two basic questions formulated at the beginning of this introduction, for this particular problem. The proposed ``geometric''  numerical integrator, which is based on splitting the orthogonal projection onto the tangent space, is robust to the (ubiquitous) presence of small singular values in the approximations. This numerically important robustness property relies on a geometric property: The low-rank manifold  is a ruled manifold (like a hyperboloid). It contains flat subspaces along which one can pass between any two points on the manifold, and the numerical integrator does just that. In this way the high curvature of the low-rank manifold at matrices with small singular values does not become harmful.
Finally, we address the nontrivial extension to tensors of various formats (Tucker tensors, tensor trains, hierarchical tensors), which is of interest in time-dependent problems with several spatial dimensions.

Section~\ref{sec:qd} on tensor and tensor network approximations in {\it quantum dynamics} combines the worlds of the previous two sections and connects them with recent developments in computational quantum physics. The reduction of the time-dependent many-particle Schr\"odinger equation to a low-rank tensor manifold by the Dirac--Frenkel time-dependent variational principle uses a tangent-space projection that is both orthogonal and symplectic. It results in a (non-canonical) Hamiltonian differential equation on the tensor manifold that can be discretized in time by the projector-splitting integrator of Section~\ref{sec:dlr}, which is robust to small singular values and preserves both the norm and the energy of the wavefunction.


\section{Hamiltonian systems of ordinary differential equations}
\label{sec:ham-ode}


\subsection{Hamiltonian systems}
\label{subsec:ham}
Differential equations of the form (with $\dot{\phantom x}=\d/\d t$)
\begin{equation} \label{ham-ode}
 \dot p = - \nabla_q H(p,q),\quad \dot q = + \nabla_p H(p,q)
\end{equation}
are fundamental to many branches of physics. The real-valued Hamilton function~$H$, defined on a domain of $\real^{d+d}$ (the phase space), represents the total energy and $q(t)\in\real^d$ and $p(t)\in\real^d$ represent the  positions and momenta, respectively, of a conservative system at time~$t$. The total energy is conserved:
$$
H(p(t),q(t)) = H(p(0),q(0))
$$
along every solution $(p(t),q(t))$ of the Hamiltonian differential equations.

\begin{figure}[t]
\centering
 \begin{picture}(0,0)
  \epsfig{file=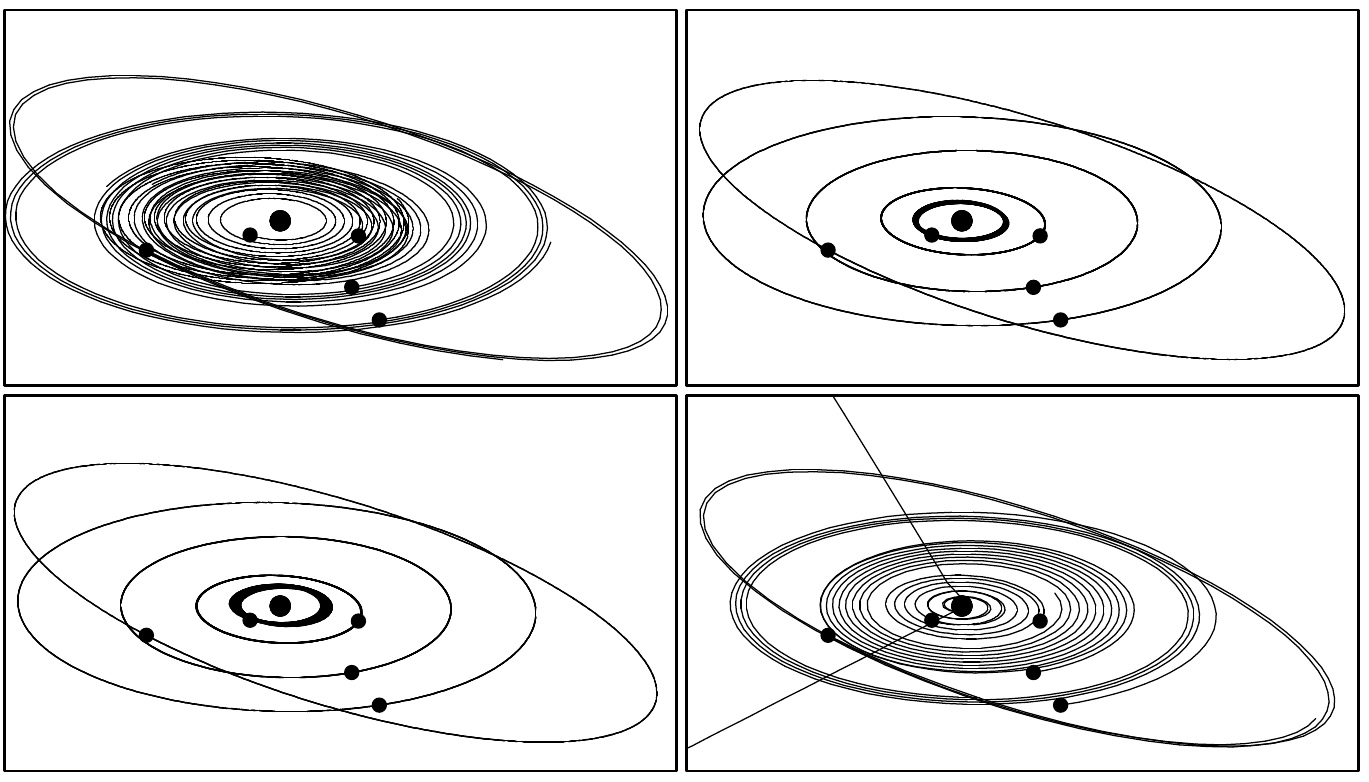}
 \end{picture}%
\begin{picture}(395.4,224.6)(  -1.4, 374.2)
  \GGGput[0,0]( 76.00,523.36){J}
  \GGGput[0,0](109.25,527.93){S}
  \GGGput[0,0](107.32,513.12){U}
  \GGGput[0,0](115.15,503.72){N}
  \GGGput[0,0]( 48.18,523.84){P}
  \GGGput[  100,  100](175.50,574.43){$h=15$}
  \GGGput[    0,  100]( 16.86,590.69){explicit Euler,  $\alpha=0,\beta =0$}
  \GGGput[0,0](272.34,523.36){J}
  \GGGput[0,0](305.58,527.93){S}
  \GGGput[0,0](303.65,513.12){U}
  \GGGput[0,0](311.48,503.72){N}
  \GGGput[0,0](244.51,523.84){P}
  \GGGput[  100,  100](371.83,574.43){$h=150$}
  \GGGput[    0,  100](213.20,590.69){symplectic Euler,  $\alpha=0,\beta =1$}
  \GGGput[0,0]( 76.00,412.42){J}
  \GGGput[0,0](109.25,416.88){S}
  \GGGput[0,0](107.32,402.18){U}
  \GGGput[0,0](115.15,392.79){N}
  \GGGput[0,0]( 48.18,412.90){P}
  \GGGput[  100,  100](175.50,463.49){$h=150$}
  \GGGput[    0,  100]( 16.86,479.63){symplectic Euler,  $\alpha=1,\beta =0$}
  \GGGput[0,0](272.34,412.42){J}
  \GGGput[0,0](305.58,416.88){S}
  \GGGput[0,0](303.65,402.18){U}
  \GGGput[0,0](311.48,392.79){N}
  \GGGput[0,0](244.51,412.90){P}
  \GGGput[  100,  100](371.83,463.49){$h=15$}
  \GGGput[    0,  100](213.20,479.63){implicit Euler,  $\alpha=1,\beta =1$}
 \end{picture}
\vspace{-2mm}
\caption{Numerical simulation of the outer solar system.}
\label{fig:solar}
\end{figure}

{\it Numerical example:} We consider four variants of the Euler method, which for a given (small) step size $h>0$ compute approximations $p_n\approx p(nh)$, $q_n\approx q(nh)$ via 
$$
p_{n+1} = p_n - h \nabla_q H(p_{n+\alpha},q_{n+\beta}) , \quad\ q_{n+1} = q_n + h \nabla_p H(p_{n+\alpha},q_{n+\beta}),
$$
with $\alpha,\beta\in\{0,1\}$. For $\alpha=\beta=0$ this is the explicit Euler method, for $\alpha=\beta=1$ it is the implicit Euler method. The partitioned methods with $\alpha\ne\beta$ are known as the {\it symplectic Euler methods}. All four methods are of order~$r=1$, that is, the error after one step of the method is $\bigo(h^{r+1})$ with $r=1$.

\begin{figure}[t]
\centering
 \begin{picture}(0,0)
  \epsfig{file=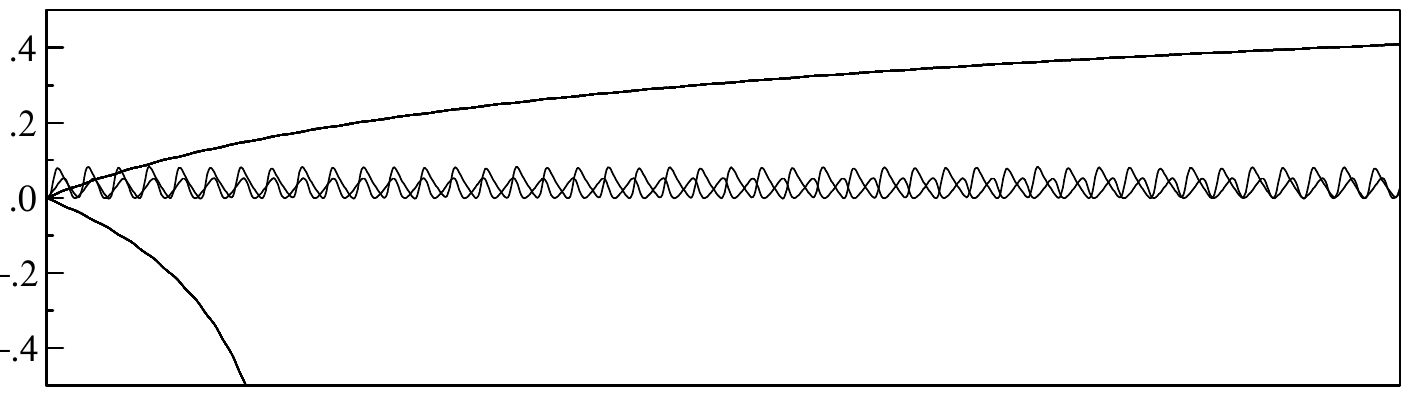}
 \end{picture}%
\begin{picture}(406.9,113.7)( -12.9, 485.2)
  \GGGput[  0,    0]( 11.08,582.01){relative error of the Hamiltonian}
  \GGGput[   50,  100](254.75,538.77){symplectic Euler,  $h = 200$}
  \GGGput[  100,    0](367.73,566.84){explicit Euler,  $h = 2$}
  \GGGput[0,0]( 48.18,511.67){implicit Euler,  $h = 2$}
 \end{picture}
\vspace{-2mm}
\caption{Relative error of the Hamiltonian on the interval
$0\le t \le 200\,000$.}
\label{fig:ham-solar}
\end{figure}

We apply these methods to the outer solar system, which is an
$N$-body problem with Hamiltonian
\[
H(p,q) = \frac 12 \sum_{i=0}^N \frac 1{m_i} \, | p^i|^2 - G
\sum_{i=1}^N \sum_{j=0}^{i-1} \frac{m_im_j}{|q^i -q^j | },
\]
where $p=(p^0,\ldots ,p^N)$, $q=(q^0,\ldots ,q^N)$ and $|\cdot|$ denotes the Euclidean norm,
and the constants are taken from \cite[Section~I.2.4]{hairer06gni}.
The positions $q^i\in\real^3$ and momenta $p^i\in\real^3$ are those of the
sun and the five outer planets (including Pluto). Figure~\ref{fig:solar} shows the
numerical solution obtained by the four versions of the Euler method on a time
interval of $200\,000$ earth days.
For the explicit Euler method the planets spiral outwards, for the implicit
Euler method they spiral inwards, fall into the sun and finally are ejected. Both
symplectic Euler methods show a qualitatively correct behaviour, even with
a step size (in days) that is much larger than the one used for the explicit and implicit
Euler methods. Figure~\ref{fig:ham-solar} shows the relative error of the
Hamiltonian, $\bigl( H(p_n,q_n) - H(p_0,q_0)\bigr)\big/ |H(p_0,q_0)|$, along the numerical
solution of the four
versions of Euler's method on the time interval $0\le nh \le 200\,000$.
Whereas the size of the error increases for the explicit and implicit Euler methods, it remains
bounded and small, of a size proportional to the step size $h$, for both symplectic Euler methods.

\subsection{Symplecticity of the flow and symplectic integrators}
\label{subsec:symp}

The time-$t$ {\it flow} of a differential equation $\dot y=f(y)$ is the map
$\phi_t$ that associates with an initial value $y_0$ at time $0$ the solution
value at time $t$: $\phi_t(y_0)=y(t)$. Consider now the Hamiltonian system~\eqref{ham-ode}
or equivalently, for $y=(p,q)$,
$$
\dot y = J^{-1}\nabla H(y) \with J=\begin{pmatrix}0 & I\\-I & 0
\end{pmatrix}.
$$
The flow $\phi_t$ of the Hamiltonian system is {\it symplectic} (or {\it canonical}), that is, the
derivative matrix $D\phi_t$ with respect to the initial value satisfies
$$
D\phi_t(y)^{\top} \,J\, D\phi_t(y) =J
$$
for all $y$ and $t$ for which $\phi_t(y)$ exists. This  quadratic
relation is formally similar to orthogonality, with $J$ in place of the
identity matrix $I$, but it is related to the preservation of areas rather
than lengths in phase space. 

There is also a local converse: If the flow of some differential equation is symplectic, then there exists locally a Hamilton function for which the corresponding Hamiltonian system coincides with this differential equation.

A numerical one-step method $y_{n+1}=\Phi_h(y_n)$ (with step size $h$) is called symplectic if
the numerical flow $\Phi_h$ is a symplectic map:
$$
D\Phi_h(y)^{\top} \,J\, D\Phi_h(y) =J.
$$
Such methods exist: the ``symplectic Euler methods'' of the previous subsection are indeed symplectic. This was first noted, or considered noteworthy, in an unpublished report by de Vogelaere \cite{vogelaere56moi}. The symplecticity can be readily verified by direct calculation or by observing that the symplectic Euler methods are symplectic maps with the $h$-scaled Hamilton function taken as the generating function  of  a canonical transformation in Hamilton and Jacobi's theory. More than 25 years later, Ruth~\cite{ruth83aci}  and Feng Kang~\cite{feng85ods,feng86dsf} independently constructed higher-order symplectic integrators using generating functions of Hamilton--Jacobi theory. These symplectic methods require, however, higher derivatives of the Hamilton function. Symplectic integrators began to find widespread interest in numerical analysis when in 1988 Lasagni,  Sanz-Serna and Suris \cite{lasagni88crm,sanz-serna88rsf,suris88otc} 
independently characterized symplectic Runge--Kutta methods by a quadratic relation of the method coefficients. This relation was already known to be satisfied by the class of Gauss--Butcher methods (the order-preserving extension of Gaussian quadrature formulae to differential equations), which include methods of arbitrary order. Like the Euler methods, Runge--Kutta methods only require evaluations of the vector field, but no higher derivatives.

The standard integrator of molecular
dynamics, introduced to the field by Verlet in 1967 \cite{verlet67ceo} and used ever since, is also symplectic.
For a Hamiltonian $H(p,q)=\frac12 p^\top M^{-1}p + V(q)$ with a symmetric positive definite mass matrix $M$, the method is explicit and given by the formulas  
\begin{eqnarray*}	
 p_{n+1/2} &=& p_n -\sfrac h2 \nabla V(q_n)\\
 q_{n+1} &=& q_n + h M^{-1}p_{n+1/2}\\
 p_{n+1} &=& p_{n+1/2} -\sfrac h2 \nabla V(q_{n+1}).
\end{eqnarray*} 
Such a method was also formulated by the astronomer
St\"ormer in 1907, and can even be traced back to  {Newton's {\it
Principia} from 1687, where it was used as a theoretical tool} in the proof of the preservation of angular momentum in the two-body problem (Kepler's second law), which is preserved by this method (cf.~\cite{wanner10kna}).
The above method is referred to as the {St\"ormer--Verlet} method, 
Verlet method or leapfrog method in different communities.
The symplecticity of this method can be understood in various ways by relating the method to classes of methods that have
proven useful in a variety of applications (cf.~\cite{hairer03gni}): as a {\it composition method} (it is a composition of the two symplectic Euler methods with half step size), as a {\it splitting method} (it solves in an alternating way the Hamiltonian differential equations corresponding to the kinetic energy $\frac12 p^\top M^{-1}p$ and the potential energy $V(q)$), •and as a {\it variational integrator\/}: it  minimizes the discrete action functional that results from approximating the action integral
$$
\int_{t_0}^{t_N} \!\! L(q(t),\dot q(t))\, \d t \with
L(q,\dot q)=\sfrac12 \dot q^\top\! M\dot q - V(q)
$$
by the trapezoidal rule 
and using piecewise linear approximation to $q(t)$.
The St\"ormer--Verlet method can thus be interpreted as resulting from a discretization of the Hamilton variational principle. Such an interpretation can in fact be given for every symplectic method. Conversely, symplectic methods can be {\it constructed} by minimizing a discrete action integral. In particular, approximating the action integral by a quadrature formula and the positions $q(t)$ by a piecewise polynomial leads to a symplectic partitioned Runge--Kutta method. With Gauss quadrature, this gives a reinterpretation of the Gauss--Butcher methods (cf.\,\cite{suris90hmo,marsden01dma,hairer06gni}).

\subsection{Backward error analysis} 
\label{subsec:bea}
The numerical example of Section~\ref{subsec:ham}, and many more examples in the literature, show that symplectic integrators behave much better over long times than their non-symplectic counterparts. How can this be explained, or put differently: How does the geometry lead to favourable dynamics? There is a caveat: As was noted early on \cite{gladman91sif,calvo92vsf}, all the benefits of symplectic integrators are lost when they are used with variable step sizes as obtained by standard step size control. So it is not just about preserving symplecticity. 

Much insight into this question is obtained from the viewpoint of {\it backward analysis}, where the result of one step of a numerical method for a differential equation $\dot y=f(y)$ is interpreted as the solution to a modified differential equation (or more precisely formal solution, having the same expansion in powers of the step size~$h$)
$$
\dot {\widetilde y} = f(\widetilde y) + hf_1(\widetilde y)+h^2f_2(\widetilde y) + h^3f_3(\widetilde y)+\ldots\ .
$$
The question then is how geometric properties of the numerical method, such as symplecticity, are reflected in the modified differential equation. It turns out that in the case of a {\it symplectic} integrator applied to a Hamiltonian differential equation, each of the perturbation terms is a {\it Hamiltonian} vector field, 
$f_j(y)= J^{-1}\nabla H_j(y)$ (at least locally, on simply connected domains). The formal construction was first given by Moser \cite{moser68loh}, where the problem of interpolating a near-identity symplectic map by a Hamiltonian flow was considered. For the important class of symplectic partitioned Runge--Kutta methods (which includes all the examples mentioned in Section~\ref{subsec:symp}), a different construction in \cite{hairer94bao}, using the theory of P-series and their associated trees, showed that the perturbation Hamiltonians $H_j$ are indeed {\it global}, defined  on the same domain on which the Hamilton function $H$ is defined and smooth.  Alternatively, this can also be shown using the explicit generating functions for symplectic partitioned Runge--Kutta methods as derived by Lasagni; see \cite[Sect.\,IX.3]{hairer06gni}. This global result is in particular important for studying the behaviour of symplectic integrators for near-integrable Hamiltonian systems, which are considered in neighbourhoods of tori. It allows us to bring the rich arsenal of Hamiltonian perturbation theory to bear on the long-time analysis of symplectic integrators.

The step from a formal theory (with the three dots at the end of the line) to rigorous estimates was taken by Benettin \& Giorgilli \cite{benettin94oth} (see also \cite{hairer97tlo,reich99bea} and \cite[Chapter IX]{hairer06gni} for related later work), who showed that in the case of an analytic vector field~$f$,  the result $y_1=\Phi_h(y_0)$ of one step of the numerical method and the time-$h$ flow $\widetilde\phi_h(y_0)$ of the corresponding modified differential equation, suitably truncated after $N\sim 1/h$ terms,  differ by a term that is exponentially small in~$1/h$:
$$
\|\Phi_h(y_0) - \widetilde\phi_h(y_0) \| \le Ch \,\e^{-c/h},
$$
uniformly for $y_0$ varying in a compact set. The constants $C$ and $c$ can be given explicitly. It turns out that $c$ is inversely proportional to a local Lipschitz constant $L$ of $f$, and hence the estimate is meaningful only under the condition $hL\ll 1$. We note that in an oscillatory Hamiltonian system, $L$ corresponds to the highest frequency in the system.

A different approach to constructing a modified Hamiltonian whose flow is exponentially close to the near-identity symplectic map is outlined by Neishtadt \cite{neishtadt84tso}, who exactly embeds the symplectic map into the flow of a non-autonomous Hamiltonian system with rapid oscillations and then uses averaging to obtain an autonomous modified Hamiltonian.

\subsection{Long-time near-conservation of energy} The above results immediately explain the observed near-preservation of the total energy by symplectic integrators used with constant step size: Over each time step, and as long as the numerical solution stays in a  fixed compact set, the Hamilton function $\widetilde H$ of the optimally truncated modified differential equation is almost conserved up to errors of size~$\bigo(h \e^{-c/h})$. On writing $\widetilde H(y_n)- \widetilde H(y_0)$ as a telescoping sum and adding up the errors, we thus obtain
$$
\widetilde H(y_n)- \widetilde H(y_0) = \bigo(\e^{-c/2h}) \for nh \le \e^{c/2h}.
$$
For a symplectic method of order $r$, the modified Hamilton function $\widetilde H$ is $\bigo(h^r)$ close to the original Hamilton function $H$, uniformly on compact sets, and so we have {\it near-conservation of energy over exponentially long times\/}:
$$
H(y_n)-  H(y_0) = \bigo(h^r) \for nh \le \e^{c/2h}.
$$
Symplecticity is, however, not necessary for good energy behaviour of a numerical method. First, the assumption can clearly be weakened to conjugate symplecticity, that is, the one-step method $y_{n+1}=\Phi_h(y_n)$ is such that $\Phi_h = \chi_h^{-1}\circ \Psi_h \circ \chi_h$  where the map $\Psi_h$ is symplectic.  But then, for some methods such as the St\"ormer--Verlet method, long-time near-conservation of energy can be proved by an argument that does not use symplecticity, but just the time-symmetry $\Phi_{-h}\circ\Phi_h={\rm id}$ of the method  \cite{hairer03gni}. That proof is similar in spirit to proving the conservation of the  energy $\frac12 p^\top M^{-1}p + V(q)=\tfrac12 \dot q^\top M\dot  q + V(q)$ for the second-order differential equation $M\ddot q +\nabla V(q)=0$ by taking the inner product with $\dot q$ and noting that there results a total  differential: $\frac{\d}{\d t} (\tfrac12 \dot q^\top M\dot  q + V(q))=0$. This kind of argument can be extended to proving long-time near-conservation of energy and momentum for symmetric multistep methods \cite{hairer04smm,hairer17*smm}, which do not preserve symplecticity and are not conjugate to a symplectic method. 
Arguments of this type are also basic in studying the long-time behaviour of numerical methods in situations where the product of the step size with the highest frequency is not very small, contrary to the condition $hL\ll 1$ above. We will encounter  such situations in Sections~\ref{sec:ham-osc} and~\ref{sec:ham-pde}.

\subsection{Integrable and near-integrable Hamiltonian systems} Symplectic integrators enjoy remarkable properties when they are used on integrable and perturbed integrable Hamiltonian systems. Their study
combines backward error analysis
and the perturbation theory of integrable
systems, a rich mathematical theory originally developed for
problems of celestial mechanics \cite{poincare92lmn,siegel71vuh,arnold97mao}.

A Hamiltonian system with the (real-analytic) Hamilton function $H(p,q)$
is called {\it integrable} if there
exists a symplectic transformation
$
(p,q) = \psi(a,\theta)
$ 
to {\it action-angle variables} $(a,\theta)$, defined for actions 
$a=(a_1,\dots,a_d)$
in some open set of $\real^d$ and for angles $\theta$ on the
$d$-dimensional torus
$
\torus^d=
\{ (\theta_1, \dots, \theta_d);\; \theta_i \in \real \hbox{ mod } 2\pi \},
$ 
such that the Hamiltonian in these variables depends only on the actions:
$$
H(p,q) = H(\psi(a,\theta)) = K(a).
$$
In the action-angle variables, the equations of motion are simply
$\dot a = 0,
\dot \theta = \omega(a)$ 
with the {\it frequencies} $\omega=(\omega_1,\dots,\omega_d)^T=\nabla_a K$.
For every $a$, the torus $\{(a,\theta): \theta\in\torus^d\}$ is thus
invariant under the flow. 
We express the actions and angles in terms of the original variables
$(p,q)$ via the inverse transform  as
$$
(a,\theta) = (I(p,q), \Theta(p,q))
$$
and note that the components of $I=(I_1,\dots,I_d)$ 
are first integrals (conserved quantities) of the integrable
system. 

The effect of a small perturbation of an integrable system
is well under control
in subsets of the phase space where the frequencies $\omega$ satisfy
Siegel's {\it diophantine condition}:
$$
|k\cdot\omega| \ge \gamma |k|^{-\nu} \quad\hbox{ for all }\ k\in\integer^d, \: k\ne 0,
$$
for some positive constants $\gamma$ and $\nu$, with 
$|k|=\sum_i |k_i|$. For $\nu>d-1$, almost all
frequencies (in the sense of Lebesgue measure) satisfy
this non-resonance condition for some $\gamma>0$. For any choice of $\gamma$ and
$\nu$ the complementary set is,
however, open and dense in $\real^d$.

For general numerical integrators applied to integrable systems 
(or perturbations thereof) the error grows quadratically with time,
and there is a linear drift in the actions $I_i$ along the numerical solution.
Consider now a symplectic partitioned Runge--Kutta method of order $r$ (or more generally, a symplectic method that has a globally defined modified Hamilton function),  applied to the integrable system with a sufficiently small step size $h\le h_0$.
Then, there
is the following result on {\it linear error growth} and {\it long-time near-preservation of
the actions}   \cite[Sect.\,X.3]{hairer02gni}: every numerical solution $(p_n,q_n)$ 
starting with frequencies
$\omega_0=\omega(I(p_0,q_0))$ such that
$\| \omega_0 - \omega^* \| \le c |\log h|^{-\nu-1}$ for some $\omega^*\in\real^d$ that obeys the above diophantine
condition, satisfies
$$
\begin{array}{rl}
\| (p_n,q_n) - (p(t),q(t)) \| &\!\!\! \le \, C\, t \, h^r \\[1mm]
\| I(p_n,q_n) - I(p_0,q_0) \| &\!\!\! \le \, C\, h^r
\end{array}
\quad\hbox{ for }\ t=nh \le h^{-r}~.
$$
(The constants $h_0,c,C$ depend on $d,\gamma,\nu$ and on bounds of the 
Hamiltonian.)
Under stronger conditions on the initial values, 
the near-preservation
of the action variables along the numerical solution holds for 
times that are exponentially long in a negative power of the step size
\cite{hairer97tlo}.
For a Cantor set of initial values and a Cantor set of step sizes
this holds even perpetually, as the existence of invariant
tori of the numerical integrator
close to the invariant tori of the integrable system was shown by Shang \cite{shang99kto,shang00rad}.

The linear error growth persists when
the symplectic integrator is applied to a
{\it perturbed integrable system} $H(p,q)+\varepsilon G(p,q)$ with
a perturbation parameter of size $\varepsilon=\bigo (h^\alpha)$ for some positive
exponent $\alpha$. 
Perturbed integrable systems have KAM tori, i.e.,
deformations of the invariant tori of the integrable system
corresponding to diophantine frequencies $\omega$, which are
invariant under the flow of the perturbed system.
If the method is applied to such a perturbed integrable system,
then the numerical method has almost-invariant tori over
exponentially long times \cite{hairer97tlo}.
For a Cantor set of non-resonant
step sizes there are even truly invariant tori
on which the numerical one-step map reduces to
rotation by $h\omega$ in suitable coordinates
\cite[Sect.\,X.6.2]{hairer02gni}.

In a very different line of research, one asks for integrable discretizations of integrable systems; see the monumental treatise by Suris~\cite{suris03tpo}.

\subsection{Hamiltonian systems on manifolds} \label{subsec:ham-mf}
In a more general setting, a Hamiltonian system is considered on a symplectic manifold, which is a manifold $\calM$ with a closed, non-degenerate alternating two-form $\omega$, called the symplectic form. Given a smooth Hamilton function $H:\calM \to \R$, the corresponding Hamiltonian differential equation is to find $u:[0,T]\to \calM$ such that
$$
\omega_{u(t)}(\dot u(t), v) = \d H(u(t))[v]   \quad\ \hbox{ for all } v \in T_{u(t)}\calM,
$$
where $T_u\calM$ denotes the tangent space at $u$ of $\calM$, for a given initial value $u(0)=u_0\in\calM$. On inserting $v=\dot u(t)$ it is seen that the total energy $H(u(t))$ is constant in time. We write again $u(t)=\phi_t(u_0)$ to indicate the dependence on the initial value. The flow map $\phi_t$ is {\it symplectic} in the sense that the symplectic form $\omega$ is preserved along the flow: for all $t$ and $u_0$ where $\phi_t(u_0)$ exists,
$$
\omega_{\phi_t(u_0)} (\d\phi_t(u_0)[\xi],\d\phi_t(u_0)[\eta]) = \omega_{u_0} (\xi,\eta)   \quad\hbox{for all } \xi,\eta \in T_{u_0}\calM; \quad\hbox{ or } \ \phi_t^*\omega=\omega.
$$
Contrary to the canonical Hamiltonian systems considered before, no general prescription is known how to construct a symplectic numerical integrator for a general Hamiltonian system on a general symplectic manifold. 

However, for the important class of Hamiltonian systems with holonomic constraints, there exist symplectic extensions of the St\"ormer--Verlet method \cite{andersen83rav,leimkuhler94sio} and of higher-order partitioned Runge--Kutta methods \cite{jay96spr}. Here the symplectic manifold $\calM$ is the submanifold of $\R^{2d}$ given by constraints $g(q)=0$, which  constrain only the positions, together with the implied constraints for the momenta, $Dg(q)\nabla_p H(p,q)=0$.

Apart from holonomic mechanical systems, there exist specially tailored symplectic integrators  for particular problem classes of non-canonical Hamiltonian systems. These are often splitting methods, as for example, for rigid body dynamics \cite{dullweber97ssm,benettin01acs}, for Gaussian wavepackets in quantum dynamics \cite{faou06api}, and for post-New\-to\-nian equations in general relativity \cite{lubich10sio}. 

\pagebreak[3]

\section{Hamiltonian systems with multiple time scales}
\label{sec:ham-osc}

\subsection{Oscillatory Hamiltonian systems}
\label{subsec:ham-osc}

The numerical experiment by Fermi, Pasta and Ulam in 1955, which showed unexpected recurrent behaviour instead of relaxation to equipartition of energy in a chain of weakly nonlinearly coupled particles,
has spurred a wealth of research in both mathematics and physics; see, e.g.,  \cite{gallavotti08tfp,berman05tfp,ford92tfp,weissert97tgo}. Even today, there are only few rigorous mathematical results for large particle numbers in the FPU problem over long times  \cite{bambusi06omi,hairer12ote}, and rigorous theory is lagging behind the insight obtained from carefully conducted numerical experiments \cite{benettin13tfp}.

\begin{figure}[!t]
\centering
 \begin{picture}(0,0)
  \epsfig{file=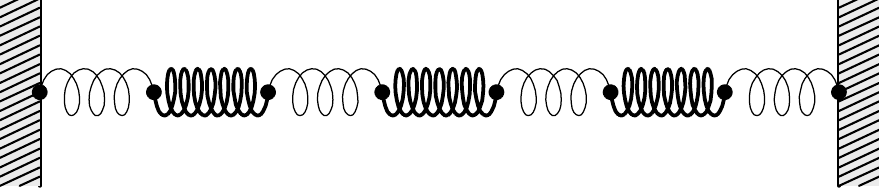}
 \end{picture}%
\begin{picture}(253.2, 54.0)(   0.0, 543.5)
  \GGGput[   50,    0]( 60.83,552.87){\small stiff}
  \GGGput[   50,    0]( 60.83,543.47){\small harmonic}
  \GGGput[   50,    0](159.48,552.87){\small soft}
  \GGGput[   50,    0](159.48,543.47){\small anharmonic}
 \end{picture}
\vspace{-2mm}
\caption{Chain of alternating stiff harmonic and soft anharmonic springs.}
\label{fig:chain}
\end{figure}

Here we consider a related class of oscillatory Hamiltonian systems for which the long-time behaviour is by now quite well understood analytically both for the continuous problem and its numerical discretizations, and which show interesting behaviour on several time scales.
The considered multiscale Hamiltonian systems couple high-frequency harmonic oscillators with a Hamiltonian of slow motion. An illustrative example of such a Hamiltonian is
provided by a Fermi--Pasta--Ulam type system of point masses
interconnected by stiff harmonic and soft anharmonic springs, as shown
in Figure~\ref{fig:chain}; see \cite{galgani92otp} and \cite[Section~I.5]{hairer06gni}.
The general setting is as follows:
For positions
$q = (q_{0},q_{1},\ldots ,q_{m})$ and momenta $p = (p_{0},p_{1},\ldots ,p_{m})$ with $p_{j},q_{j}\in\real^{d_{j}}$, let the Hamilton function be given by
$$
H(p ,q ) =  H_{\omega}(p ,q ) + H_{\rm slow} (p ,q ) ,
$$
where the oscillatory and slow-motion energies are given by
$$
H_{\omega}(p ,q )  =  \sum_{j=1}^m \sfrac 12\Bigl( |p_{j}|^2 +
 \omega_{j}^2 \,|q_{j}|^2\Bigr), \qquad
H_{\rm slow}(p ,q ) = \sfrac 12  |p_{0}|^2 + U(q )
$$
with high frequencies
$$
\omega_{j} \ge \eps^{-1}, \qquad
0<\eps\ll 1 .
$$
The coupling potential $U$ is assumed smooth with derivatives bounded independently of the small parameter $\eps$.
On eliminating the momenta $p_j=\dot q_j$, the Hamilton equations become the system of second-order differential equations
$$
\ddot q_j + \omega_j^2 q_j =
 -\nabla_j U(q ) , \qquad j=0,\dots,m,
$$
where $\nabla_j$ denotes the gradient with respect to $q_j$, and where we set $\omega_{0}=0$. We are interested in the behaviour of the system for initial values with an oscillatory energy that is bounded independently of $\eps$:
$$ 
H_{\omega}(p(0) ,q(0) ) \le {\rm Const.}
$$
This system shows different behaviour on different time scales:
\begin{enumerate}[(i)]
\item
almost-harmonic motion of the fast variables $(p_j,q_j)$ $(j\ne 0)$  on time scale $\eps$;
\item
motion of the slow variables $(p_0,q_0)$  on the time scale $\eps^0$; 
\item
energy exchange between
the harmonic oscillators with the same frequency on the time scale $\eps^{-1}$;
\item
energy exchange between
the harmonic oscillators corresponding to frequencies in 1:2 or 1:3 resonance on the time scale $\eps^{-2}$ or $\eps^{-3}$, respectively;
\item
near-preservation of the $j$th oscillatory energy $E_j=\frac 12 ( |p_{j}|^2 +
 \omega_{j}^2 \,|q_{j}|^2)$ for a non-resonant frequency $\omega_j$ beyond the time scale $\eps^{-N}$ for arbitrary $N$; and
\item
near-preservation of the total oscillatory energy $H_\omega$
over intervals that are beyond the time scale $\eps^{-N}$ for arbitrary $N$, uniformly for $\omega_j\ge\eps^{-1}$ without any non-resonance condition (but depending on the number $m$ of different frequencies). Hence, there is nearly no energy exchange between the slow and the fast subsystem for  very long times irrespective of resonances, almost-resonances or non-resonances among the high frequencies.
\end{enumerate}
The long-time results (iii)--(vi) require in addition that $q_0$ stays in a compact set, which is ensured if the potential $U(q)\to+\infty$ as $|q|\to\infty$. These results can be proved by two alternative techniques: 
\begin{itemize}
\item[(H)]
 using canonical coordinate transformations of Hamiltonian perturbation theory to a normal form; or
\item[(F)]
using modulated Fourier expansions in time.
\end{itemize}
The latter technique was developed by the authors of the present paper (in part together with David Cohen) and will be outlined in the next subsection. 

Motivated by the problem of relaxation times in statistical mechanics, item (v) was first shown using (H) by Benettin, Galgani \& Giorgilli~\cite{benettin87roh} for the single-frequency case, even over times exponentially long in $\eps^{-1}$ for a real-analytic Hamilton function, and in \cite{benettin89roh} for the multi-frequency case over times exponentially long in some negative power of $\eps$ that depends on the diophantine non-resonance condition; it was subsequently shown using (F) in \cite{cohen03mfe} over exponentially long times in $\eps^{-1}$ for the single frequency-case, and in  \cite{cohen05nec} over times  $\eps^{-N}$ with $N$ depending on the chosen non-resonance condition on the frequencies.

Item (vi) was first shown using (F) in \cite{gauckler13esi} and subsequently using (H) in \cite{bambusi13nfa}, where the result was extended to exponentially long time scales. 

The relationship between the two techniques of proof, (H) and (F), is not clear at present. The proofs look very different in the basic arguments, in the geometric content and in the technical details, yet lead to very similar results about the long-time behaviour of the continuous problem.

\subsection{Modulated Fourier expansion}
\label{subsec:mfe}
Modulated Fourier expansions in time have proven useful in the long-time analysis of  differential equations where the nonlinearity appears as a perturbation to a primarily linear problem (as laid out in the programme of \cite{fermi55} cited in the introduction). This encompasses important classes of Hamiltonian ordinary and partial differential equations. The approach can be successfully used for the analysis of the continuous problems as well as for their numerical discretizations, as is amply shown in the corresponding references in this and the next section. In particular for the analysis of numerical methods, it offers the advantage that it does not require nonlinear coordinate transforms. Instead, it embeds the original system in a high-dimensional system of modulation equations that has a Lagrangian / Hamiltonian structure with invariance properties.
 In addition to the use of modulated Fourier expansions as an analytical technique, they have been used also as a numerical approximation method in \cite[Chapter XIII]{hairer02gni} and \cite{cohen04diss,condon09oho,condon10oso,faou14aps,bao14aua,zhao17uam}.

We now describe the basic steps how, for the problem of the previous subsection, a simple ansatz for the solution over a short time interval leads to long-time near-conservation results for the oscillatory energies $E_j=\frac 12 ( |p_{j}|^2 +\omega_{j}^2 \,|q_{j}|^2)$.
We approximate the solution $q_j$ of the second-order differential equation of the previous section as a {\it modulated Fourier expansion},
$$
 q_j(t) \approx \sum_{\bfk} z_j^\bfk(t)\, \e^{\iu(\bfk\cdot{\bfomega})t} \qquad \hbox{for short times $0\le t \le 1$},
$$
with modulation functions $z_j^\bfk$, all derivatives of which are required to be bounded independently of $\eps$. The sum is taken over a finite set of
multi-indices
$\bfk=(k_1,\dots,k_m)\in \Z^m$, and 
$\bfk\cdot{\bfomega}=\sum k_j\omega_j$. The slowly changing modulation functions are multiplied with the highly oscillatory exponentials $\e^{\iu(\bfk\cdot{\bfomega})t}=\prod_{j=1}^m \bigl(\e^{\iu\omega_j t}\bigr)^{k_j}$, which are products of solutions to the linear equations $\ddot x_j + \omega_j^2 x_j = 0$. Such products can be expected to be introduced into the solution $q_j$ by the nonlinearity. 

Similar multiscale expansions have appeared on various occasions in the literature. The distinguishing feature here is that such a short-time expansion is used to derive long-time properties of the Hamiltonian system.
%

\subsubsection{Modulation system and non-resonance condition}
When we insert this ansatz into the differential equation and collect the coefficients to the same exponential $\e^{\iu(\bfk\cdot{\bfomega})t}$, we obtain the infinite system of modulation equations for $\bfz=(z_j^\bfk)$
$$
    {(\omega_{j}^{2}-(\bfk\cdot\bfomega)^{2})}\,{z_{j}^{\bfk}}+ 2\iu(\bfk\cdot\bfomega){\dot z}_{j}^{\bfk} +{\ddot z}_{j}^{\bfk}    = -\frac{\partial\,\calU}{\partial z_{j}^{-\bfk}}\,(\bfz).
$$
 The left-hand side results from the linear part $\ddot q_j + \omega_j^2 q_j$ of the differential equation. The right-hand side results from the nonlinearity and turns out to have a gradient structure 
with the modulation potential
\begin{eqnarray*}
&&
\calU (\bfz ) = 
U(z^\bfzero ) 
+ \sum_{\ell\ge 1}\
\sum_{\bfk^1+\ldots +\bfk^\ell =\bfzero}\frac 1 {\ell!}\, U^{(\ell)}(z^\bfzero )
\bigl[ z^{\bfk^1},\ldots ,z^{\bfk^\ell} \bigr],
\end{eqnarray*}
with $0\ne\bfk^i\in \Z^m$. The sum is suitably truncated, say to $\ell\le N$, $|\bfk^i|\le N$.
%

The infinite modulation system is truncated and can be 
solved approximately (up to a defect~$\eps^N$) 
for modulation functions $z_j^\bfk$ with derivatives bounded independently of $\eps$ under a non-resonance condition that ensures that $\omega_{j}^{2}-(\bfk\cdot\bfomega)^{2}$ is the dominating coefficient on the left-hand side of the modulation equation,  except when $k$ is plus or minus the $j$th unit vector. For example, we can assume, as in \cite{cohen05nec}, that there exists $c>0$ such that
$$
|\bfk\cdot\bfomega| \ge c\,\eps^{-1/2} \for \bfk\in \Z^m \with 0<|\bfk|\le 2N,
$$
where $|\bfk| = \sum_{j=1}^m |k_j|$. Under such a non-resonance condition one can construct and appropriately bound the modulation functions $z_j^\bfk$, and the modulated Fourier expansion 
is then an $\bigo(\eps^N)$ approximation to the solution over a short time interval  $t = \bigo(1)$.
%

\subsubsection{Lagrangian structure and invariants of the modulation system.}
With the multiplied  functions $y_j^\bfk(t)=z_j^\bfk(t)\e^{\iu(\bfk\cdot{\bfomega})t}$ that appear as summands in the modulated Fourier expansion, the modulation equations take the simpler form
$$
\ddot y_j^\bfk +\omega_j^2 y_j^\bfk = -\frac{\partial\,\calU}{\partial y_{j}^{-\bfk}}\,(\bfy).
$$
The modulation potential $\calU$ has the important invariance property
$$
\calU(S_{\ell}(\theta)\bfy)=\calU(\bfy) \for
S_{\ell}(\theta)\bfy   = (\e^{\iu k_{\ell} \theta}y^{\bfk}_j)_{j,\bfk} \with \theta\in \R,
$$
as is directly seen from the sum over ${\bfk^1+\ldots +\bfk^m =\bfzero}$ in the definition of $\calU$. We have thus embedded the original differential equations in a system of modulation equations that are Lagrange equations for the Lagrange function
$$
\calL(\bfy, \dot\bfy) = \sfrac12 \sum_{j,\bfk} \dot y_j^{-\bfk} \dot y_j^\bfk -  \sfrac12 \sum_{j,\bfk} \omega_j^2 y_j^{-\bfk}  y_j^\bfk-\calU(\bfy),
$$
which is invariant under the action of the group $\{ S_{\ell}(\theta): \theta\in\R \}$. We are thus in the realm of Emmy Noether's theorem from 1918, which states that invariance under a continuous group action (a geometric property) yields the existence of conserved quantities of the motion (a dynamic property).
By Noether's theorem, the modulation equations thus conserve
$$
\calE_{\ell}(\bfy,\dot\bfy) 
  =  -\iu\sum_{j}\sum_{\bfk} k_{\ell}\omega_{\ell} \, y_j^{-\bfk} \, \dot y_j^\bfk.
$$
Since the modulation equations are solved only up to a defect $\bigo(\eps^N)$ in the construction of the modulated Fourier expansion, the functions $\calE_{\ell}$ are almost-conserved quantities with $\bigo(\eps^{N+1})$ deviations over intervals of length $\bigo(1)$. They turn out to be $\bigo(\eps)$ close to the oscillatory energies $E_{\ell}$. By patching together many short time intervals, the drift in the almost-invariants $\calE_{\ell}$ is controlled to remain bounded by $Ct\eps^{N+1} \le C\eps$  over long times $t\le \eps^{-N}$, and hence also the deviation in the oscillatory energies $E_{\ell}$ is  only $\bigo(\eps)$ over such long times. We 
thus obtain long-time near-conservation of the oscillatory energies $E_{\ell}$.

\subsection{Long-time results for numerical integrators}
\label{subsec:ham-osc-num}
Modulated Fourier expansions were first developed in \cite{hairer00lec} and further in \cite[Chapter XIII]{hairer02gni} to understand the observed long-term near-conservation of energy by some numerical methods for step sizes for which the smallness condition $hL\ll 1$ of the backward error analysis of Section~\ref{subsec:bea} is not fulfilled.
For the numerical solution of the differential equation of Section~\ref{subsec:ham-osc}, we are interested in using  numerical integrators that allow large step sizes $h$ such that $h/\eps\ge c_0>0$. In this situation, the one-step map of a numerical integrator is no longer a near-identity map, as was the case in Section~\ref{sec:ham-ode}.

For a class of time-symmetric {\it trigonometric integrators}, which are exact  for the uncoupled harmonic oscillator equations $\ddot x_j + \omega_j^2 x_j=0$ and reduce to the St\"ormer--Verlet method for $\omega_j=0$, the
following results are proved for step sizes $h$ that satisfy a {\it numerical non-resonance condition\/}:
$$
\text{$h\omega_j$ is bounded away (by $\sqrt{h}$) from a multiple of $\pi$.}
$$
Under just this condition it is shown in \cite{cohen15lta}, using modulated Fourier expansions, that the slow energy $H_{\rm slow}$ is nearly preserved along the numerical solution for very long times $t\le h^{-N}$ for arbitrary $N\ge 1$ provided the total energy remains bounded along the numerical solution. If in addition,
$$
\begin{aligned}[l]
&\text{sums of $\pm h\omega_j$ with at most $N+1$ terms are} \\[-1mm]
 &\text{bounded away from non-zero multiples of $2\pi$,}
\end{aligned}
$$
then also the total and oscillatory energies $H$ and $H_\omega$ are 
nearly preserved along the numerical solution for  $t\le h^{-N}$ for the symplectic methods among the considered symmetric trigonometric integrators. Modified total and oscillatory energies are nearly preserved by the non-symplectic methods in this class.
These results yield the numerical version of property (vi) above. 
A numerical version of property (v) was shown in \cite{cohen05nec}. The single-frequency case was previously studied in \cite{hairer00lec}. For the  {\it St\"ormer--Verlet method}, which  can be interpreted as a trigonometric integrator with modified frequencies, related long-time results are given in \cite{hairer00ecb,cohen15lta}.

The numerical version of the energy transfer of property (iii) was studied in \cite[Section XIII.4]{hairer02gni} and in \cite{cohen05nec,mclachlan14mti}. Getting the energy transfer qualitatively correct by the numerical method turns out to put more restrictions on the choice of methods than long-time energy conservation.

%

While we concentrated here on long-time results, it should be mentioned that fixed-time convergence results of numerical methods for the multiscale problem as $h\to 0$ and $\eps\to 0$ with $h/\eps\ge c_0>0$ also pose many challenges; see, e.g., \cite{garcia-archilla99lmf,hochbruck99agm,grimm06eao,buchholz17ctg} for systems with constant high frequencies and also \cite{lubich14nif,hairer16lta} for systems with state-dependent high frequencies, where near-preservation of adiabatic invariants is essential. We also refer to \cite[Chapters XIII and XIV]{hairer06gni} and to the review \cite{cohen06nif}.

%
%
\section{Hamiltonian partial differential equations}
\label{sec:ham-pde}
There is a vast literature on the long-time behaviour of nonlinear wave equations, nonlinear Schr\"odinger equations and other Hamiltonian partial differential equations; see, e.g., the monographs \cite{kuksin93nii,bourgain99gso,craig00pdp,kuksin00aoh,kappeler03kdv,
grebert14tdn} where infinite-dimensional versions of Hamiltonian perturbation theory are developed. Here we consider a few analytical results that have recently been transfered also to numerical discretizations.

\subsection{Long-time regularity preservation}

We consider the nonlinear wave equation (or nonlinear Klein--Gordon equation)
\[
\partial_t^2 u = \partial_x^2 u - \rho u + g(u), \qquad u=u(x,t)\in\mathbb{R}
\]
with $2\pi$-periodic boundary condition in one space dimension, a positive mass parameter $\rho$ and a smooth nonlinearity $g=G'$ with $g(0)=g'(0)=0$.
This equation is a Hamiltonian partial differential equation
$
\partial_t v = - \nabla_u H(u,v)$, $\partial_t u = \nabla_v H(u,v)
$
(where $v=\partial_t u$)
with Hamilton function
\[
H(u,v) = \frac{1}{2\pi} \int_{-\pi}^\pi \biggl( \sfrac12 \Bigl( v^2 + (\partial_x u)^2 + \rho u^2\Bigr) - G(u) \biggr) \,\d x
\]
on the Sobolev space $H^1$ of $2\pi$-periodic functions.

Written in terms of the Fourier coefficients $u_j$ of $u(x,t) = \sum_{j\in\integer} u_j(t) \e^{\iu jx}$, the nonlinear wave equation takes the form of the oscillatory second-order differential equation of Section~\ref{subsec:ham-osc}, but the system is now infinite-dimensional:
\[
\ddot{u}_j + \omega_j^2 u_j = \mathcal{F}_j g(u), \qquad j\in\integer,
\]
where $\mathcal{F}_j$ gives the $j$th Fourier coefficient and 
$
\omega_j = \sqrt{j^2 + \rho} 
$
are the frequencies.

The following result is proved, using infinite-dimensional Hamiltonian perturbation theory (Birkhoff normal forms), by Bambusi~\cite{bambusi03bnf}, for arbitrary $N\ge 1$: Under a non-resonance condition on the frequencies $\omega_j$, which is satisfied for almost all values of the parameter $\rho$, 
and for initial data $(u^0,v^0)$ that are $\eps$-small in a Sobolev space $H^{s+1}\times H^s$ with sufficiently large $s=s(N)$, 
the harmonic energies $E_j=\frac 12 ( |\dot{u}_{j}|^2 + \omega_{j}^2 \,|u_{j}|^2)$ are nearly preserved over the time scale $t\le \eps^{-N}$, and so is the $H^{s+1}\times H^s$ norm of the solution $(u(t),v(t))$. 

An alternative proof using modulated Fourier expansions was given in \cite{cohen08lta} with the view towards transfering the result to numerical discretizations with trigonometric integrators as done in \cite{cohen08coe}, for which in addition also a numerical non-resonance condition is required.
%

Related long-time near-conservation results are proved for other classes of Hamiltonian differential equations, in particular for nonlinear Schr\"odinger equations with a resonance-removing convolution potential,
in \cite{bourgain96coa,bambusi06bnf,grebert07bnf} using Birkhoff normal forms and in \cite{gauckler10nse} using modulated Fourier expansions. These results are transfered to numerical discretization by Fourier collocation in space and splitting methods in time in  \cite{faou10bna,faou10bnb,gauckler10sif}. 

For small initial data where only one pair of Fourier modes is excited (that is, $E_j(0)=0$ for $|j|\ne 1$, which would yield a plane wave solution in the linear wave equation), higher Fourier modes become excited in the above nonlinearly perturbed wave equation to yield, within short time, mode energies $E_j(t)$ that decay geometrically with $|j|$ and then remain almost constant for very long times. Using modulated Fourier expansions, this is proved in \cite{gauckler12mes} for the continuous problem  and extended to numerical discretizations in
\cite{gauckler17mes}.

\subsection{Long-time near-conservation of energy for numerical discretizations}
In Section~\ref{sec:ham-ode} we have seen that symplectic integrators approximately preserve the energy of a Hamiltonian ordinary differential equation over long times. It is then not far-fetched to expect the same for full (space and time) discretizations of Hamiltonian partial differential equations. However, the standard backward error analysis argument does {\it not} carry over from ODEs to PDEs:

-- Backward error analysis requires that the product of the time step size, now denoted $\dt$, and the local Lipschitz constant $L$ be  small. After space discretization with a meshwidth $\dx$, we have $L$ proportional to $1/\dx$ in the case of  wave equations and proportional to $1/\dx^2$ in the case of  Schr\"odinger equations. The condition $\dt \, L \ll 1$ then requires unrealistically small time steps as compared with the familiar Courant--Friedrichs--Lewy (CFL) condition $\dt \,L \le \mathrm{Const}$.

-- Even when a symplectic integrator is used with a very small time step size such that $\dt \, L \ll 1$, the numerical method will nearly preserve the Hamilton function of the spatially discretized system, not that of the PDE. The two energies are close to each other only as long as the numerical solution is sufficiently regular, which usually cannot be guaranteed {\it a priori}.

The above hurdles are overcome in \cite{cohen08coe} for the nonlinearly perturbed wave equation of the previous subsection discretized by Fourier collocation in space and symplectic trigonometric integrators in time. Here, high regularity of the numerical solution and near-conservation of energy are proved simultaneously using modulated Fourier expansions. In \cite{gauckler10sif}, this technique and the energy conservation results are taken further to a class of nonlinear Schr\"odinger equations with a resonance-removing convolution potential (in arbitrary space dimension) discretized by Fourier collocation in space and a splitting method in time.

Long-time near-conservation of energy for symplectic splitting methods applied to the nonlinear Schr\"odinger equation in one space dimension (without a resonance-removing convolution potential) is shown in \cite{faou11hio,faou12gni} with a backward error analysis adapted to partial differential equations and, under weaker step size restrictions,  in \cite{gauckler16nlt} with modulated Fourier expansions. In contrast to the aforementioned results, these results are not restricted to initial values in higher-order Sobolev spaces.

Apart from these results in the weakly nonlinear regime, the basic question of long-time approximate conservation of the total energy under numerical discretization remains wide open in view of the difficulties addressed above.


\subsection{Orbital stability results for nonlinear Schr\"odinger equations and their numerical discretizations}
We consider the cubic nonlinear Schr\"odinger equation 
$$
i \partial_t u = -\Delta u + \kappa |u|^2 u
$$
near special solutions (ground states and plane waves) and describe results about orbital stability that have been obtained for the continuous problem and for numerical discretizations.

{\it Ground state.}
The cubic nonlinear Schr\"odinger equation  on the real line in 
the focusing case $\kappa=-1$ admits  the solitary wave solution $u(x,t)= \e^{\iu \lambda t} \eta(x)$ with $\eta(x)=\tfrac1{\sqrt{2}}\mathrm{sech}(x/2)$ with an appropriate real parameter $\lambda$. It is known from Weinstein \cite{weinstein85mso} that this solution  is orbitally stable in the sense that for a small $H^1$-perturbation of the initial data, the exact solution remains close to the orbit of $\eta$ for all times. For the case restricted to symmetric initial conditions, orbital stability of a full discretization with a finite-difference space discretization and a splitting method as time discretization is shown by Bambusi, Faou \& Gr\'ebert \cite{bambusi13eas}.

{\it Plane waves.}
We now consider the cubic nonlinear Schr\"odinger equation on a torus $\mathbb{T}^d=(\real/2\pi\Z)^d$ of arbitrary dimension $d\ge 1$ with real $\kappa$. It admits the plane wave solution $u(x,t)=\rho \e^{\iu (m\cdot x)}\e^{-\iu (|m|^2+\kappa|\rho|^2) t}$ for arbitrary $\rho>0$ and $m\in\Z^d$. It is shown in \cite{faou13ssp} that for almost all $\rho>0$,
 plane wave solutions are orbitally stable for long times $t\le \eps^{-N}$ under $\eps$-perturbations in higher-order Sobolev spaces $H^s$, $s=s(N)\gg 1$, if $1+2\kappa|\rho|^2>0$ (which is the condition for linear stability). This result is given with two different proofs, one based on Birkhoff normal forms and the other one on modulated Fourier expansions.  The latter technique is used in \cite{faou14pws} to transfer the result to numerical discretization using Fourier collocation in space and a splitting method for time discretization. The long-time orbital stability under smooth perturbations is in contrast to the instability under rough perturbations shown in \cite{hani14lti}.

\section{Dynamical low-rank approximation}
\label{sec:dlr}
Low-rank approximation of too large matrices and tensors is a fundamental approach to data compression and model reduction in a wide range of application areas. Given a matrix $A\in\R^{m\times n}$, the best rank-$r$ approximation to $A$ with respect to the distance given by the Frobenius norm (that is, the Euclidean norm of the vector of entries of a matrix) is known to be obtained by a truncated singular value decomposition:
$
A \approx \sum_{i=1}^r \sigma_i u_i v_i^\top,
$
where $\sigma_1,\dots,\sigma_r$ are the $r$ largest singular values of $A$, and $u_i\in\R^m$ and $v_i\in\R^n$ are the corresponding left and right singular vectors, which form an orthonormal basis of the range and corange, respectively, of the best approximation. Hence, only $r$ vectors of both length $m$ and $n$ need to be stored. If $r\ll\min(m,n)$, then the requirements for storing and handling the data are significantly reduced.

When $A(t)\in\R^{m\times n}$, $0\le t \le T$, is a time-dependent family of large matrices, computing the best rank-$r$ approximation would require singular value decompositions of $A(t)$ for every time instance $t$ of interest, which is often not computationally feasible. Moreover, when $A(t)$ is not given explicitly but is the unknown solution to a matrix differential equation $\dot A(t)=F(t,A(t))$, then computing the best rank-$r$ approximation would require to first solve the differential equation on $\R^{m\times n}$, which may not be  feasible for large $m$ and $n$, and then to compute the singular value decompositions at all times of interest, which may again not be feasible.

\subsection{Dynamical low-rank approximation of matrices}

An alternative --- and often computationally feasible --- approach can be traced back to Dirac \cite{dirac30noe} in a particular context of quantum dynamics (see also the next section). Its abstract version can be viewed as a nonlinear Galerkin method on the tangent bundle of an approximation manifold $\calM$ and reads as follows: Consider a differential equation $\dot A(t)=F(t,A(t))$  in a (finite- or infinite-dimensional) Hilbert space $\calH$, and let $\calM$ be a submanifold of $\calH$. An approximation $Y(t)\in\calM$ to a solution $A(t)$ (for $0\le t \le T$) is determined by choosing the time derivative $\dot Y(t)$ as the orthogonal projection of the vector field $F(t,Y(t))$ to the tangent space $T_{Y(t)}\calM$ at $Y(t)\in\M$:
\begin{equation}\label{dlr}
\dot Y(t) = P_{Y(t)} F(t,Y(t)),
\end{equation}
where $P_Y$ denotes the orthogonal projection onto the tangent space at $Y\in\calM$. Equation~\eqref{dlr} is a differential equation on the approximation manifold $\calM$, which is complemented with an initial approximation $Y(0)\in\calM$ to $A(0)\in\calH$. When $\calM$ is a flat space, then this is the standard Galerkin method, which is a basic approximation method for the spatial discretization of partial differential equations.  When $\calM$ is not flat, then the tangent space projection $P_Y$ depends on $Y$, and \eqref{dlr} is a nonlinear differential equation even if $F$ is linear.

For the dynamical low-rank approximation of time-dependent matrices, \eqref{dlr} is used with $\calM$ chosen as the manifold of rank-$r$ matrices in the space $\calH=\R^{m\times n}$ equipped with the Frobenius inner product (the Euclidean inner product of the matrix entries). This approach was first proposed and studied in \cite{koch07dlr}. The rank-$r$ matrices are represented in (non-unique) factorized form as
$$
Y=USV^\top,
$$
where $U\in\R^{m\times r}$ and $V\in\R^{n\times r}$ have orthonormal columns and $S\in\R^{r\times r}$ is an invertible matrix. The intermediate small matrix $S$ is not assumed diagonal, but it has the same non-zero singular values as $Y\in\M$. Differential equations for the factors $U,S,V$ can be derived from \eqref{dlr} (uniquely under the gauge conditions $U^\top\dot U=0$ and $V^\top\dot V=0$). They contain the inverse of $S$ as a factor on the right-hand side. It is a typical situation that $S$ has small singular values, because in order to obtain accurate approximability, the discarded singular values need to be small, and then the smallest retained singular values are usually not much larger. Small singular values complicate the analysis of the approximation properties of the dynamical low-rank approximation \eqref{dlr}, for a geometric reason: the curvature of the rank-$r$ manifold $\M$ at $Y\in\M$ (measured as the local Lipschitz constant of the projection map $Y\mapsto P_Y$) is proportional to the inverse of the smallest singular value of $Y$. It seems obvious that high curvature of the approximation manifold can impair the approximation properties of \eqref{dlr}, and for a general manifold this is indeed the case. Nevertheless, for the  manifold $\M$ of rank-$r$ matrices there are numerical and theoretical results in \cite{koch07dlr} that show good approximation properties also in the presence of arbitrarily small singular values.

\subsection{Projector-splitting integrator} \label{subsec:psi}

The numerical solution of the differential equations for $U,S,V$ encounters difficulties with standard time integration methods (such as explicit or implicit Runge--Kutta methods) when $S$ has small singular values, since the inverse of $S$ appears as a factor on the right-hand side of the system of differential equations. 

A numerical integration method for these differential equations with remarkable properties is given in \cite{lubich14aps}. It is based on splitting the tangent space projection, which at $Y=USV^\top$ is an alternating sum of three subprojections:
$$
P_Y Z = ZVV^\top -  UU^\top Z  VV^\top + UU^\top Z.
$$
Starting from a factorization $Y_n=U_nS_nV_n^\top$ at time $t_n$, the corresponding splitting integrator updates the factorization of the rank-$r$ approximation to $Y_{n+1}=U_{n+1}S_{n+1}V_{n+1}^\top$ at time $t_{n+1}$. It alternates between solving (approximately if need be) matrix differential equations of dimensions $m\times r$ (for $US$), $r\times r$ (for $S$), $n\times r$ (for $VS^\top$) and doing orthogonal decompositions of matrices of these dimensions. The inverse of $S$ does not show up in these computations.

The projector-splitting  integrator has a surprising exactness property: if the given matrix $A(t)$ is already of rank $r$ for all $t$, then the integrator reproduces $A(t)$ exactly \cite{lubich14aps}. More importantly, the projector-splitting  integrator is robust to the presence of small singular values: it admits convergent error bounds that are independent of the singular values \cite{kieri16ddl}. The proof uses the above exactness property and a geometric peculiarity: in each substep of the algorithm, the approximation moves along a flat subspace of the manifold $\M$ of rank-$r$ matrices. In this way, the high curvature due to small singular values does no harm.

\subsection{Dynamical low-rank approximation of tensors}
The dynamical low-rank approximation and the projector-splitting integrator have been extended from matrices to tensors $A(t)\in\R^{n_1\times\dots\times n_d}$ such that the favourable approximation and robustness properties are retained; see \cite{koch10dta,lubich13dab,arnold14ota,lubich15tio,lubich17tio}. The dynamical low-rank approximation can be done in various tensor formats that allow for a notion of rank, such as Tucker tensors, tensor trains, hierarchical tensors, and general tensor tree networks; see \cite{hackbusch12tsa,uschmajew13tga} for these concepts and for some of their geometric properties.

\section{Quantum dynamics}
\label{sec:qd}
\subsection{The time-dependent variational approximation principle}
The time-dependent Schr\"odinger equation for the $N$-particle wavefunction $\psi=\psi(x_1,\dots,x_N,t)$,
$$
\iu \partial_t \psi = H\psi,
$$
posed as an evolution equation on the complex Hilbert space $\calH=L^2((\R^{3})^N,\C)$ with a self-adjoint Hamiltonian operator $H$, is not accessible to direct numerical treatment in the case of several, let alone many particles. ``One must therefore resort to approximate methods'', as Dirac \cite{dirac30noe} noted already in the early days of quantum mechanics. For a particular approximation scheme, which is nowadays known as the time-dependent Hartree--Fock method, he used the tangent space projection \eqref{dlr} for the Schr\"odinger equation. Only later was this recognized as a general approximation approach, which is now known as the (Dirac--Frenkel) time-dependent variational principle in the physical and chemical literature: Given a submanifold $\M$ of $\calH$, an approximation $u(t)\in\M$ to the wavefunction $\psi(\cdot,t)\in\calH$ is determined by the condition that
$$
\dot u \text{ is chosen as that $w\in T_u\M$ for which } \| \iu w -  Hu \| \text{ is minimal.}
$$
This is precisely \eqref{dlr} in the context of the Schr\"odinger equation: $\dot u = P_u \frac1\iu Hu$. If we assume that the approximation manifold $\M$ is such that for all $u\in\M$,
$$
T_u\M \text{ is a complex vector space},
$$
(so that with $v\in T_u\M$, also $\iu v\in T_u\M$), then the orthogonal projection $P_u$ turns out to be also a {\it symplectic} projection with respect to the canonical symplectic two-form on $\calH$ given by $\omega(\xi,\eta)= 2 \Imag \langle \xi,\eta \rangle$ for $\xi,\eta\in\calH$, and $\M$ is a symplectic manifold.
With the Hamilton function $H(u)=\langle u, Hu \rangle$, the differential equation for $u$ can then be rewritten as
$$
\omega(\dot u,v) = \d H(u)[v] \qquad\text{for all }v\in T_u\M,
$$
which is a Hamiltonian system on the symplectic manifold $\M$; cf.~Section~\ref{subsec:ham-mf}. The total energy $H(u)$ is therefore conserved along solutions, and the flow is symplectic on $\M$. The norm is conserved if $\M$ contains rays, i.e., with $u\in\M$ also $\alpha u\in\M$ for all ${\alpha>0}$.
We refer to the books \cite{kramer81got,lubich08fqt} for geometric, dynamic and approximation aspects of the time-dependent variational approximation principle.

\subsection{Tensor and tensor network approximations} 
In an approach that builds on the time-honoured idea of separation of variables, the multi-configuration time-dependent Hartree method (MCTDH) \cite{meyer90tmc,meyer09mqd} uses the time-dependent variational principle to determine an approximation to the multivariate wavefunction that is a linear combination of products of univariate functions:
$$
u(x_1,\dots,x_N,t) = \sum_{i_1=1}^{r_1} \dots \sum_{i_N=1}^{r_N} c_{i_1,\dots,i_N}(t) \,\varphi_{i_1}^{(1)}(x_1,t)\dots \varphi_{i_N}^{(N)}(x_N,t).
$$
The time-dependent variational principle yields a coupled system of ordinary differential equations for the coefficient tensor $\bigl(c_{i_1,\dots,i_N}(t)\bigr)$ of full multilinear rank and low-dimensional nonlinear Schr\"odinger equations for the single-particle functions $\varphi_{i_n}^{(n)}(x_n,t)$, which are assumed orthonormal for each $n=1,\dots,N$. Well-posedness and regularity for this nonlinear system of evolution equations is studied in \cite{koch07rot}, and an asymptotic error analysis of the MCTDH approximation for growing ranks $r_n$ is given in \cite{conte10aea}. 

The projector-splitting integrator of Section~\ref{subsec:psi} is extended to MCTDH in \cite{lubich15tii}. The nonlinear MCTDH equations are thus split into
a chain of linear single-particle differential equations, alternating with orthogonal matrix decompositions. The integrator conserves the $L^2$ norm and the total energy and, as is proved in \cite{lubich17tio}, it is robust to the presence of small singular values in matricizations of the coefficient tensor.

In the last decade, tensor network approximations, and in particular matrix product states, have  increasingly come into use for the description of strongly interacting quantum many-body  systems; see, e.g., \cite{verstraete08mps,cirac09rat,szalay15tpm}. Matrix product states (known as tensor trains in the mathematical literature \cite{oseledets11ttd}) approximate the wavefunction by
$$
u(x_1,\dots,x_N,t) =G_1(x_1,t)\cdot\ldots\cdot G_N(x_N,t)
$$
with matrices $G_n(x_n,t)$ of compatible (low) dimensions. This approach can be viewed as a non-commutative separation of variables. Its memory requirements grow only linearly with the number of particles $N$, which makes the approach computationally attractive for many-body systems. The approximability of the wavefunction or derived quantities by this approach is a different issue, with some excellent computational results but hardly any rigorous mathematical theory so far.

For the numerical integration of the equations of motion that result from the time-dependent variational approximation principle, the projector-splitting integrator has recently been extended to matrix product states in \cite{lubich15tio,haegeman16ute}, with favourable properties like the MCTDH integrator. The important robustness to the presence of small singular values is proved in \cite{kieri16ddl}, again using the property that the integrator moves along flat subspaces within the tensor manifold.

\subsection*{Acknowledgement}
We thank Bal\'azs Kov\'acs, Frank Loose, Hanna Walach, and Gerhard Wanner for helpful comments.


\newcommand{\etalchar}[1]{$^{#1}$}

\end{document}